\documentclass[11pt]{article}
\usepackage{amssymb}
\newtheorem{lemma}{Lemma}[section]

\newtheorem{theorem}[lemma]{Theorem}
\newtheorem{proposition}[lemma]{Proposition}
\newtheorem{corollary}[lemma]{Corollary}

\newtheorem{remark}[lemma]{Remark}

\def\endproof{\hfill$\Box$}

\def\endproof{\hfill$\Box$}

\title{Spherically orderable groups\footnote{The work was carried out in the framework of the
State Contract of the Sobolev Institute of Mathematics, Project
No.~FWNF-2022-0012.}}
\author{S.V. Sudoplatov}
\date{}

\begin{document}

\maketitle
\begin{abstract}
We introduce and study the class of spherically ordered groups.
The notions of spherically ordered groups and their spectra of
spherical orderability are introduced. Values of these spectra are
found for a series of natural groups.
\end{abstract}


\bigskip
\section{Introduction}

Well known linearly ordered groups and circularly ordered groups
are both deeply investigated and described \cite{Fuchs, KoKo} and
admit various generalizations and modifications for partial, left
and right  orderings \cite{Fuchs, KoMe}, betweenness and
separation groups \cite{Shepperd1, Shepperd2, Shepperd3} and
semigroups \cite{Gilder}.

We continue to study $n$-spherical orders introducing and
investigating $n$-sphe\-ri\-cal\-ly ordered groups which
generalize linearly and cyclically ordered groups, with $n=2$ and
$n=3$, respectively. The notions of spherically ordered groups and
their spectra of spherical orderability are introduced. Values of
these spectra are found for a series of natural groups.

\section{Spherically orderable groups and their spherical spectra}

Recall the following generalization of linear and circular orders.

\medskip
{\bf Definition} \cite{aafot, dso, mcso}. An $n$-ary relation
$K_n$ is called {\em $n$-ball}, or {\em $n$-spherical}, or {\em
$n$-circular} order relation, for $n\geq 2$, if it satisfies the
following conditions:

(nso1)\footnote{Here we admit cyclic permutations iff $n$ is odd. In view of this circumstance the item (nso1) is corrected with respect to the previous version.} for any even permutation $\sigma$ on $\{1,2,\ldots,n\}$, $$\forall x_1,\ldots,x_n \left(K_n(x_1,x_2,\ldots,x_n)\to K_n\left(x_{\sigma(1)},x_{\sigma(2)}\ldots,x_{\sigma(n)}\right)\right);$$\\

(nso2) $\forall x_1,\ldots,x_n
\biggm((K_n(x_1,\ldots,x_i,\ldots,x_j,\ldots,x_n)\land$
$$\land K_n(x_1,\ldots,x_j,\ldots,x_i,\ldots,x_n))\leftrightarrow\bigvee\limits_{1\leq k<l\leq n} x_k= x_l\biggm)$$ for any $1\leq i<j\leq n$;\\

(nso3) $\forall x_1,\ldots,x_n\Biggm(K_n(x_1,\ldots,x_n)\to$
$$\to\forall t\left(\bigvee\limits_{i=1}^nK_n(x_1,\ldots,x_{i-1},t,x_{i+1},\ldots,x_n)\right)\Biggm);$$\\

(nso4) $\forall x_1,\ldots,x_n
(K_n(x_1,\ldots,x_i,\ldots,x_j,\ldots,x_n)\lor$
$$\lor K_n(x_1,\ldots,x_j,\ldots,x_i,\ldots,x_n)),\,\, 1\leq i<j\leq
 n.$$

A structure $\mathcal{M}$ provided with a $n$-spherical order is
called {\em $n$-spherically ordered}.

\medskip
The axioms above produce all possible linear orders $K_2$ and
circular orders $K_3$. Here (nso2) gives the reflexivity: $\forall
x K_2(x,x)$, and the antisymmetry: $$\forall
x_1,x_2(K_2(x_1,x_2)\wedge K_2(x_2,x_1)\to x_1\approx x_2),$$
(nso2) and (nso3) give the transitivity: $$\forall
x_1,x_2,x_3(K_2(x_1,x_2)\wedge K_2(x_2,x_3)\to K_2(x_1,x_3)),$$
and the axiom (nso4) gives the linearity: $\forall
x_1,x_2(K_2(x_1,x_2)\vee K_2(x_2,x_1)).$ For the transitivity it
suffices to take pairwise distinct elements $a,b,c$ with $(a,b)\in
K_2$ and $(b,c)\in K_2$. By (nso3) we have $(a,c)\in K_2$ or
$(c,b)\in K_2$, and $(b,a)\in K_2$ or $(a,c)\in K_2$. But the
cases $(c,b)\in K_2$ and $(b,a)\in K_2$ are impossible in view of
(nso2), implying the required $(a,c)\in K_2$.

Similarly to linear orders any $n$-spherical order $K_n$ on a set
$Z$ has the dual one consisting of all $n$-tuples in $Z^n\setminus
K_n$ united with the set of all $n$-tuples in $Z^n$ with some
repeated coordinates. We denote this dual order by
$\overline{K_n}$.

\medskip
{\bf Definition.} A group $G$ is called {\em $n$-spherically
ordered}, or {\em $n$-$s$-ordered}, if $G$ is provided with a
$n$-spherical order $K_n$ such that for any $(x_1,\ldots,x_n)\in
K_n$ and any $y\in G$ the tuples $(x_1y,\ldots,x_ny)$ and
$(yx_1,\ldots,yx_n)$ belong to $K_n$.

A group $G$ is called {\em $n$-spherically orderable}, or {\em
$n$-$s$-orderable}, if $G$ has a $n$-spherically ordered
expansion. A group $G$ is called {\em spherically orderable} if it
is $n$-spherically orderable for some $n$.

For a group $G$ we define its {\em spectrum ${\rm Sp}_{\rm so}$ of
spherical orderability}, or {\em spherical spectrum}, as follows:
$$
{\rm Sp}_{\rm so}(G)=\{n\in\omega\setminus\{0,1\}\mid G\mbox{ is
}n\mbox{-spherically orderable}\}.
$$

A group $G$ is called {\em totally spherically orderable}, or {\em
totally $s$-orderable}, if $G$ has maximal spectrum of spherical
orderability, i.e. ${\rm Sp}_{\rm so}(G)=\omega\setminus\{0,1\}$.

A group $G$ is called {\em almost totally spherically orderable},
or {\em almost totally $s$-orderable}, if  ${\rm Sp}_{\rm so}(G)$
is a cofinite subset of $\omega$.

A group $G$ is ({\em almost}) {\em not $s$-orderable in any way}
if ${\rm Sp}_{\rm so}(G)$ is empty (respectively, finite).

\medskip
The notions above for the spherical orderability and its spectra
can be naturally spread for an arbitrary structure. Besides, the
spherical orderability admits similar variations of orderability
as for linear (bi-)orderability such as left-orderability and
right-orderability \cite{KoMe, Mul}.

A natural {\bf problem} arises on description of spherical spectra
for groups and related structures.

\medskip
By the definition a group $G$ is linearly ordered iff $G$ is
$2$-spherically ordered, and $G$ is cyclically ordered iff $G$ is
$3$-spherically ordered. Here $2\in{\rm Sp}_{\rm so}(G)$ and
$3\in{\rm Sp}_{\rm so}(G)$, respectively.

Again by the definition a group $G$ is spherically orderable iff
${\rm Sp}_{\rm so}(G)\ne\emptyset$.

\begin{remark}\label{rem_spher0}\rm
If a group $G$ is $n$-spherically orderable then each subgroup of
$G$ is $n$-spherically orderable, too, since any restriction
$\mathcal{M}$ of a $n$-spherically ordered structure
$\mathcal{N}$, with a $n$-spherical order $K_n$, is again
$n$-spherically ordered, with the $n$-spherical order $K_n\cap
M^n$.
\end{remark}

In view of Remark \ref{rem_spher0} we have the following {\em
Monotonicity property} for the spectrum of spherical orderability:

\begin{proposition}\label{pr_sp_mon}
For any groups $G_1, G_2$ if $G_1\leqslant G_2$ then ${\rm
Sp}_{\rm so}(G_1)\supseteq{\rm Sp}_{\rm so}(G_2)$.
\end{proposition}

Proposition \ref{pr_sp_mon} immediately implies:

\begin{corollary}\label{cor_tot}
If $G$ is an {\rm (}almost{\rm )} totally $s$-orderable group then
any its subgroup is also {\rm (}almost{\rm )} totally
$s$-orderable.
\end{corollary}

\begin{corollary}\label{cor_nordaw}
If $G$ is {\rm (}almost{\rm )} not $s$-orderable in any way then
any its supergroup is also {\rm (}almost{\rm )} not $s$-orderable
in any way.
\end{corollary}

Proposition \ref{pr_sp_mon} and Corollary \ref{cor_tot} can be
also deduced from the following criterion of $n$-spherical
orderability:

\begin{theorem}\label{th_sp_cr}
A group $G$ is $n$-spherically ordered by a $n$-spherical order
$K_n$ iff for any $n$-tuple $(a_1,\ldots,a_n)\in K_n$ with
pairwise distinct coordinates and for any $b\in G$ the tuples
$(a_1b,\ldots,a_nb)$ and $(ba_1,\ldots,ba_n)$ are even
permutations, i.e. are not odd permutations of some tuples in
$K_n$.
\end{theorem}

Proof. By the definition of $n$-spherical order $K_n$ it consists
of all $n$-tuples with some repeated coordinates and of all even
permutations of given $n$-tuples in $K_n$ such that for any
$n$-tuple, in the universe, with pairwise distinct coordinate
either this tuple or all its odd permutations belong to $K_n$. Now
if $G$ is $n$-spherically ordered by $K_n$ it is forbidden to
include into $K_n$ its odd permutations of forms
$(a_1b,\ldots,a_nb)$ and $(ba_1,\ldots,ba_n)$, and conversely
permitted to include its even permutations of forms
$(a_1b,\ldots,a_nb)$ and $(ba_1,\ldots,ba_n)$. \endproof

\begin{remark}\label{rem_sp_perm}\rm
Theorem \ref{th_sp_cr} shows that the only obstacle for a
$n$-spherical order $K_n$ on a group $G$ to produce this group to
be $n$-spherical ordered is the possibility for the
multiplications $(a_1b,\ldots,a_nb)$ and $(ba_1,\ldots,ba_n)$ to
produce an odd permutation for a tuple $(a'_1,\ldots,a'_n)$ in
$K_n$ with pairwise distinct coordinates. In particular, elements
of stabilizers for $n$-element sets $\{a_1,\ldots,a_n\}\subseteq
G$ can not produce odd permutations of tuples $(a_1,\ldots,a_n)\in
K_n$. Here $b\ne e$ and therefore the maps $a_i\mapsto a_ib$ and
$a_i\mapsto ba_i$ do not have fixed elements.

In general, $K_n$ is uniquely defined by its subrelation $K^0_n$
consisting of all $n$-tuples in $K_n$ with pairwise distinct
coordinates, since $K_n=K^0_n\,\dot{\cup}\,K^1_n$, where $K^1_n$
consists of all tuples in $G^n$ with repeated coordinates. The
subrelation $K^0_n$ produces an algebra $\mathcal{K}^0_n$ with
$\displaystyle\frac{n!}{2}$ unary operations of even permutations
forming the alternating group $A_n$, and unary operations $r_b$
and $l_b$, for $b\in G$, carrying out the maps $a_i\mapsto a_ib$
and $a_i\mapsto ba_i$, respectively. Since all these operations
are invertible, with $(r_b)^{-1}=r_{b^{-1}}$ and
$(l_b)^{-1}=l_{b^{-1}}$, and form a {\em derivative} group
$S=S(G,n)$ with the identical even permutation which is equal both
to $r_e$ and $l_e$, we obtain an $S$-act $\mathcal{K}^0_n=\langle
K^0_n, s\rangle_{s\in S}$, for $|G|\geq n$. Here actions $r_b$ and
$l_b$ can be even permutations if, for instance, $|G|=n$, which
implies $S=A_n$, They may not belong to $A_n$, if $b\ne e$ and $G$
is torsion-free. And they coincide if $b$ belongs to the center of
$G$. The group $S$ is generated by its two subgroups $A_n$ and
$S'=\{r_b,l_b\mid b\in G\}$ and elements in $f\in A_n$ and in
$s\in S'$ commute: $fs=sf$. Besides, $r_b=l_b$ and
$\alpha_{bb'}=\alpha_b \alpha_{b'}$ for any $\alpha\in\{r,l\}$ and
$b,b'\in G$, if $G$ is commutative. Since $A_n$ is commutative iff
$n\leq 3$, the algebra $\mathcal{K}^0_n$ has the commutative
derivative group $S$ iff $n\leq 3$ and $G$ is commutative.

Clearly, if $G$ is a group $n$-spherically ordered by an order
$K_n$ then $G$ is $n$-spherically ordered by the dual order
$\overline{K_n}=(G^n\setminus K_n)\cup K^1_n$, too.
\end{remark}

In view of Theorem \ref{th_sp_cr} and Remark \ref{rem_sp_perm} we
have the following:

\begin{corollary}\label{cor_lin_cr}
A group $G$ is linearly ordered by a $2$-spherical order $K_2$ iff
for any pair $(a_1,a_2)\in K_2$ with $a_1\ne a_2$ and for any
$b\in G$ the pairs $(a_1b,a_2b)$ and $(ba_1,ba_2)$ are not
transpositions of some pairs in $K_2$.
\end{corollary}

Since all odd permutations for the group $S_3$ are transpositions
we additionally have:

\begin{corollary}\label{cor_circ_cr}
A group $G$ is circularly ordered by a circular order $K_3$ iff
for any triple $(a_1,a_2,a_3)\in K^0_3$ and for any $b\in G$ the
triples $(a_1b,a_2b,a_3b)$ and $(ba_1,ba_2,ba_3)$ are not
transpositions of some triples in $K^0_3$.
\end{corollary}

\begin{remark}\label{rem_sp_perm2}\rm
Each orbit $O$ of the group $S$ on $K^0_n$ connects all even
permutations of tuples in $O$ and possibly tuples based on
distinct $n$-element sets if some $r_b$ or $l_b$ are not
permutations of these sets.

The number $r$ of these orbits $O$ is called the {\em rank of
generation} of $K^0_n$ with respect to $S$ and denoted by ${\rm
rk}_S(K^0_n)$. This rank is said to be the rank of generation for
$K_n$, denoted by ${\rm rk}_S(K_n)$.

If the group $G$ is abelian then  ${\rm rk}_S(K^0_n)$ is finite
iff $G$ is finite, since for tuples $(x_1,\ldots,x_n)$ and
$(y_1,\ldots,y_n)$ based on distinct $n$-element sets the link
$r_b=l_b$ is defined by $b=y_1-x_1$. For an infinite non-abelian
group one may put defining relations connecting $n$-tuples into
finitely many orbits via chains of tuples. Indeed, the connection
of tuples $(x_1,\ldots,x_n)$ and $(y_1,\ldots,y_n)$ can be
organized via intermediate tuples $(z_1,\ldots,z_n)$ formed by new
defining elements and with appropriate $r_b$ and $l_b$, where the
elements $b$ are composed with some new defining elements. Thus
there are infinite non-abelian groups with finite ${\rm
rk}_S(K^0_n)$. Moreover, finite values ${\rm rk}_S(K^0_n)$ can be
realized arbitrarily in $\omega\setminus\{0\}$ using defining
relations connecting step-by-step orbits by $r_b$ and $l_b$ with
new defining elements $b$.
\end{remark}

In view of Remark \ref{rem_sp_perm2} for an abelian group $G$ the
operators $r_b=l_b$ on $G^3$ either fix all coordinates of a
triple or move all its coordinates, i.e. can not generate odd
permutations (transpositions). Thus Corollary \ref{cor_circ_cr}
immediately implies:

\begin{corollary}\label{cor_circ_cr2}
For any abelian group $G$, $3\in{\rm Sp}_{\rm so}(G)$.
\end{corollary}

The following construction based on Theorem \ref{th_sp_cr} shows
how the spectrum ${\rm Sp}_{\rm so}(G)$ can be reduced till the
empty one, i.e. with $G$ which is not $s$-orderable in any way.

For a required group $G$ we take a generating set $\{a_n\mid
n\in\omega\setminus\{0\}\}\cup\{b,c\}$. Now we consider the
following defining relations: $ba_1c=a_2$, $ba_2c=a_1$,
$ba_nc=a_n$ for $n\geq 3$. These relations show that any $m$-tuple
$(ba_{i_1}c,ba_{i_2}c,\ldots,ba_{i_m}c)$ consisting of pairwise
distinct elements and with some $i_j=1$ and $i_k=2$ is a
transposition of $(a_{i_1},a_{i_2},\ldots,a_{i_m})$. Therefore
$m$-spherical orders on $G$ coordinated with left and right group
actions can not be formed for any $m$. Hence we have the
following:

\begin{theorem}\label{th_sp_empty}
There exists a group $G$ such that ${\rm Sp}_{\rm
so}(G)=\emptyset$.
\end{theorem}

{\bf Definition.} Let $<$ be a strict linear order on a group $G$,
$C$ be the {\em cyclification} of $<$ consisting of all tuples
$(x,y,z)$ with $x<y<z\vee y<z<x\vee z<x<y$. A $n$-ary relation
$K^0_n$ on $G$, for $n\geq 3$, is called {\em $<$-coordinated} if
$K^0_n$ consists of all $n$-tuples $(a_0,\ldots,a_{n-1})$ with
$a_0<\ldots <a_{n-1}$, and for any $b\in G$ the tuples
$(a_0b,\ldots,a_{n-1}b)$ and $(ba_0,\ldots,ba_{n-1})$ satisfy
$C\left(a_ib,a_{(i+1)\,({\rm mod}\,n)}b,a_{(i+2)\,({\rm
mod}\,n)}b\right)$ and $C\left(ba_i,ba_{(i+1)\,({\rm
mod}\,n)},ba_{(i+2)\,({\rm mod}\,n)}\right),$ $i=0,\ldots,n-1$.

\medskip
\begin{remark}\label{ev_od}\rm By Remark \ref{rem_sp_perm} any $<$-coordinated relation $K^0_n$ is uniquely
extensible to its closure $K_n$ under even permutations of tuples
and addition of all tuples in $G^n$ with repeated coordinates. In
view of Remark \ref{rem_sp_perm} this closure satisfies the axioms
of $n$-spherical order iff even per\-mu\-ta\-tions of tuples
$(a_1,\ldots,a_n)$ in $K_n$ do not meet their odd permutations
under group actions $(a'_1b,\ldots,a'_nb)$ and
$(ba''_1,\ldots,ba''_n)$ in $K_n$, where
$(a'_1,\ldots,a'_n),(a''_1,\ldots,a''_n)\in K_n$, $b\in G$.
\end{remark}

\begin{remark}\label{rem_spher1}\rm
If $G$ is a finite group with $|G|=m$ and $m<n$ then $G$ is
$n$-spherically orderable. Indeed, in such a case the relation
$K_n$ consisting of all tuples with repeated coordinates satisfies
the axioms of $n$-spherical order implying that $K_n=K^1_n$ and
$\langle G,K_n\rangle$ is a $n$-spherically ordered group.
\end{remark}

In view of Remark \ref{rem_spher1} we have the following:

\begin{proposition}\label{pr_sp_fin}
If $G$ is a group with $|G|=m\in\omega$ then ${\rm Sp}_{\rm
so}(G)\supseteq\{n\in\omega\mid n>m\}$, in particular, $G$ is
almost totally $s$-orderable.
\end{proposition}

The following assertion is a reformulation of well-known folklore
fact that nonunit linearly ordered groups are infinite.

\begin{proposition}\label{pr_lo}
A finite group $G$ is $2$-spherically orderable iff $|G|=1$.
\end{proposition}

Proof. The one-element group $G$ is $2$-spherically orderable in
view of Remark \ref{rem_spher1}. If $|G|>1$ then $G$ contains a
nonunit element $g$ of finite order. Assuming that $G$ is
$2$-spherically orderable with the order $K_2$ we take the unit
$e$ and obtain $(e,g)\in K_2$ or $(g,e)\in K_2$. By the definition
of $2$-spherically ordered group we obtain $(g^k,g^{k+1})\in K_2$
for any $k$ or $(g^{k+1},g^{k})\in K_2$ for any $k$. For the first
case we take $k_0>0$ with $g^{k_0+1}=e$ and obtain $(g^{k_0},e)\in
K_2$, $g^{k_0}\ne e$, which contradicts the transitivity. In the
second case we also have a contradiction to the transitivity in
view of $(e,g^{k_0})\in K_2$, $g^{k_0}\ne e$. \endproof

\medskip
The arguments for Proposition \ref{pr_lo} imply the following:

\begin{corollary}\label{cor_lo1}
Any $2$-spherically orderable group $G$ is torsion-free, that is,
if $G$ contains an element of finite positive order then
$2\notin{\rm Sp}_{\rm so}(G)$.
\end{corollary}

\begin{remark}\label{rem_Z_m}\rm
Since each group $\mathbb Z_m$ is circularly orderable
\cite{Fuchs}, i.e. $3\in {\rm Sp}_{\rm so}(\mathbb Z_m)$, then
Proposition \ref{pr_sp_fin} and Co\-rol\-lary~\ref{cor_lo1} imply
that ${\rm Sp}_{\rm so}(\mathbb Z_2)={\rm Sp}_{\rm so}(\mathbb
Z_3)=\omega\setminus\{0,1,2\}$.

The groups $\mathbb Z_m$, for even $m$, are not $m$-spherically
orderable since the definition of $m$-spherically orderable group
implies that $K_m$ is closed both under even and odd permutations
of $m$ pairwise distinct elements: odd permutations are obtained
by actions of elements in $\mathbb Z_m$ on $m$-tuples of pairwise
distinct elements of $\mathbb Z_m$. It contradicts the axiom
(nso2). Thus $m\notin {\rm Sp}_{\rm so}(\mathbb Z_m)$ for each
even $m$. In particular, by Proposition \ref{pr_sp_fin} and
Co\-rol\-lary~\ref{cor_lo1}, ${\rm Sp}_{\rm so}(\mathbb
Z_4)=\omega\setminus\{0,1,2,4\}$.

The arguments above show that any group $G$ containing a subgroup
$\mathbb Z_m$, for even $m$, can not be $m$-spherically ordered by
a relation $K_m$ since, in view of (nso4), it should contain a
tuple $(i_1,\ldots,i_m)$, where
$\{i_1,\ldots,i_m\}=\{0,\ldots,m-1\}$. In particular, if $G$
contains subgroups $\mathbb Z_m$, for each even $m$, then ${\rm
Sp}_{\rm so}(G)$ is contained in the set of odd numbers.

At the same time the groups $\mathbb Z_m$, for odd $m\geq 5$, are
$m$-spherically orderable since the definition of $m$-spherically
orderable group implies even permutations of elements in $K_m$,
generated by the tuple $(0,1,\ldots,m-1)$ only. Thus $m\in {\rm
Sp}_{\rm so}(\mathbb Z_m)$ for each odd $m\geq 5$. \endproof
\end{remark}

In view of Proposition \ref{pr_sp_mon} and Remark \ref{rem_Z_m}
each element $a$ of even order $m$ in a given group $G$ implies
$m\notin{\rm Sp}_{\rm so}(G)$. Thus we have the following:

\begin{proposition}\label{pr_so_sp1}
For any group $G$, ${\rm Sp}_{\rm so}(G)$ does not contain even
numbers which are equal to orders of elements in $G$.
\end{proposition}

\begin{proposition}\label{pr_so1}
For any natural $m, n$ with $2<n<m$, and $n\!\!\not|m$ if $n$ is
even, then the group $\mathbb Z_m$ is $n$-spherically orderable.
\end{proposition}

Proof. We form the $n$-spherical order $K_n$ on the universe
$\mathbb Z_m$ adding to the set of $n$-tuples with some repeated
coordinates all even permutations of tuples $(k_1({\rm
mod}\,m),\ldots,k_n({\rm mod}\,m))$ with $k_1<\ldots<k_n$ and even
permutations of tuples $(k_1+q({\rm mod}\,m),\ldots,k_n+q({\rm
mod}\,m))$ with $q\in\mathbb Z_m$, i.e. $K_n$ is generated by its
$<$-coordinated subrelation $K^0_n$, where $<$ is the natural
strict order on $\mathbb Z_m$.

Since $n\!\!\not|m$ for even $n$, $\mathbb Z_m$ does not contain a
subgroup $\mathbb Z_n$ violating the $n$-spherical orderability as
in Remark \ref{rem_Z_m}. Thus even permutations of tuples
$(k_1+q({\rm mod}\,m),\ldots,k_n+q({\rm mod}\,m))$ are coordinated
with even permutations of tuples $(k_1,\ldots,k_n)$, i.e. they do
not produce odd per\-mu\-ta\-tions. Since these per\-mu\-ta\-tions
cover all possibilities for tuples with even permutations only,
Theorem \ref{th_sp_cr} guarantees that $\mathbb Z_m$ is
$n$-spherically orderable by $K_n$.
\endproof

\medskip
Propositions \ref{pr_sp_fin}, \ref{pr_lo}, \ref{pr_so1} and Remark
\ref{rem_Z_m} immediately imply the following description of
spherical spectra for the groups $\mathbb Z_m$:

\begin{theorem}\label{th_Z_m}
Let $m\in\omega\setminus\{0,1\}$. Then $${\rm Sp}_{\rm so}(\mathbb
Z_m)=\omega\setminus(\{0,1,2\}\cup\{n\mid n|m\mbox{ and }n\mbox{
is even}\}).$$
\end{theorem}

\begin{proposition}\label{pr_so1_tot}
The group $\mathbb Z$ is totally $s$-orderable.
\end{proposition}

Proof. Let $n\in\omega\setminus\{0,1\}$. We form the $n$-spherical
order $K_n$ on the universe $\mathbb Z$ adding to the set of
$n$-tuples with some repeated coordinates all even permutations of
tuples $(k_1,\ldots,k_n)$ with $k_1<\ldots<k_n$, i.e. generate
$K_n$ by the $<$-coordinated relation $K^0_n$ with the natural
order $<$. Clearly, for any $m\in\mathbb Z$,
$(k_1+m,\ldots,k_n+m)$ preserves the set of these tuples and
satisfies the axioms of $n$-spherical orders in view of Theorem
\ref{th_sp_cr}, as required. Thus, ${\rm Sp}_{\rm so}(\mathbb
Z)=\omega\setminus\{0,1\}$, i.e. $\mathbb Z$ is totally
$s$-orderable. \endproof

\medskip
It is known that any torsion-free abelian group can be
lexicographically ordered with respect to its generators.
Therefore the arguments for Proposition \ref{pr_so1_tot} imply the
following:

\begin{theorem}\label{th_tf_tot}
Any torsion-free abelian group is totally $s$-orderable.
\end{theorem}

\begin{remark}\label{rem_tf}\rm
The group confirming Theorem \ref{th_sp_empty} is torsion-free and
it is not $s$-orderable at all. It illustrates that the
commutativity of a group is essential in Theorem \ref{th_tf_tot}.
\end{remark}

\end{document}